\documentclass[reqno,12pt]{amsart}
\usepackage{amsfonts}
\usepackage{bbm}
\usepackage{} %leqno is the option to put formula numbers on the left side
\numberwithin{equation}{section}
\usepackage{indentfirst}
 \usepackage{color}
\usepackage{amssymb}
\usepackage{mathrsfs}
\usepackage{xy}
\xyoption{all}
\def\Ext{\mbox{\rm Ext}\,} \def\Hom{\mbox{\rm Hom}} \def\dim{\mbox{\rm dim}\,} \def\Iso{\mbox{\rm Iso}\,}\def\Heis{\mbox{\rm Heis}\,}
\def\lr#1{\langle #1\rangle}   \def\lra{\longrightarrow} \def\mod{\mbox{\rm \textbf{mod}}\,}
     \def\ra{\rightarrow}
\def\scrP{{\mathscr{P}}}
 \def\gl.{\mbox{\rm gl.}\,}
\def\ooz{\Omega}\def\vez{\varepsilon}
\def\Aut{\mbox{\rm Aut}\,} 

\theoremstyle{plain} %text of this environment is typesetted in italics
\newtheorem{theorem}{\bf Theorem}[section]
\newtheorem{lemma}[theorem]{\bf Lemma}
\newtheorem{corollary}[theorem]{\bf Corollary}
\newtheorem{proposition}[theorem]{\bf Proposition}

\theoremstyle{definition} %text of this environment is typesetted in roman letters
\newtheorem{definition}[theorem]{\bf Definition}
\newtheorem{remark}[theorem]{\bf Remark}
\newtheorem{example}[theorem]{\bf Example}

\newcommand{\bt}{\begin{theorem}}
\newcommand{\et}{\end{theorem}}
\newcommand{\bl}{\begin{lemma}}
\newcommand{\el}{\end{lemma}}
\newcommand{\bd}{\begin{definition}}
\newcommand{\ed}{\end{definition}}
\newcommand{\bc}{\begin{corollary}}
\newcommand{\ec}{\end{corollary}}
\newcommand{\bp}{\begin{proof}}
\newcommand{\ep}{\end{proof}}
\newcommand{\bx}{\begin{example}}
\newcommand{\ex}{\end{example}}
\newcommand{\br}{\begin{remark}}
\newcommand{\er}{\end{remark}}
\newcommand{\be}{\begin{equation}}
\newcommand{\ee}{\end{equation}}
\newcommand{\ba}{\begin{align}}
\newcommand{\ea}{\end{align}}
\newcommand{\bn}{\begin{enumerate}}
\newcommand{\en}{\end{enumerate}}
\newcommand{\bcs}{\begin{cases}}
\newcommand{\ecs}{\end{cases}}
\makeatletter
\renewcommand{\section}{\@startsection{section}{1}{0mm}
  {-\baselineskip}{0.5\baselineskip}{\bf\leftline}}
\makeatother

\begin{document}

\title[Bridgeland's Hall algebras and Heisenberg doubles]{Bridgeland's Hall algebras and Heisenberg doubles} %title of paper and the running head option

\author[Haicheng Zhang]{{Haicheng Zhang}} %first author's name and the running head option

%\dedicatory{Dedicated to Professor Xxx Yyy on his sixtieth birthday}

%%%%%%%%%%%%%%% footnote %%%%%%%%%%%%%%%%
%In case \subjclass[2010] command is not effective
%(or the version of amsart.cls is old), write as follows:
%\renewcommand{\thefootnote}{\fnsymbol{footnote}}
%\footnote[0]{2010\textit{ Mathematics Subject Classification}.
%Primary 00; Secondary 00.}
%
%\thanks{ %acknowledgment of support etc. if any
%$^{*}$Thanks.
%}
%%%%%%%%%%%% Authors addresses %%%%%%%%%%%%%
\address{% First Author
Institute of Mathematics, School of Mathematical Sciences, Nanjing Normal University,
 Nanjing 210023, P. R. China.\endgraf
}

\address{
Yau Mathematical Sciences Center, Tsinghua University,
 Beijing 100084, P. R. China.\endgraf
}
\email{zhanghai14@mails.tsinghua.edu.cn}

%%%%%%%%%%%%%%%%%%%%%%%%%%%%%%%%%%%%%%%%%

\maketitle

\begin{abstract}
Let $A$ be a finite dimensional hereditary algebra over a finite field, and let $m$ be a fixed integer such that $m=0$ or $m>2$. In the present paper, we first define an algebra $L_m(A)$ associated to $A$, called the $m$-periodic lattice algebra of $A$, and then prove that it is isomorphic to Bridgeland's Hall algebra $\mathcal {D}\mathcal {H}_m(A)$ of $m$-cyclic complexes over projective $A$-modules. Moreover, we show that there is an embedding of the Heisenberg double Hall algebra of $A$ into $\mathcal {D}\mathcal {H}_m(A)$.
\end{abstract}

\section{Introduction}
The Hall algebra $\mathfrak{H}(A)$ of a finite dimensional algebra $A$ over a finite field was introduced by Ringel \cite{R90a} in 1990. Ringel \cite{R90,R90a} proved that if $A$ is representation-finite and hereditary, the twisted Hall algebra $\mathfrak{H}_{v}(A)$, called the Ringel--Hall algebra, is isomorphic to the positive part of the corresponding quantized enveloping algebra.  By introducing a bialgebra structure on $\mathfrak{H}_{v}(A)$,
Green \cite{Gr95} generalized Ringel's work to an arbitrary finite dimensional hereditary algebra $A$
and showed that the composition subalgebra of $\mathfrak{H}_{v}(A)$ generated by simple $A$-modules
gives a realization of the positive part of the quantized enveloping algebra associated with $A$.

In order to give a Hall algebra realization of the entire quantized enveloping algebra, one considers defining the Hall algebras of triangulated categories (for example, \cite{Kapranov}, \cite{Toen}, \cite{Xiao}, \cite{XiaoXu}).
Kapranov \cite{Kapranov} defined an associative algebra, called the lattice algebra, for the bounded derived category of any hereditary algebra. By using the fibre products of model categories, To\"{e}n \cite{Toen} defined an associative algebra, called the derived Hall algebra, for DG-enhanced triangulated categories.
Later on, Xiao and Xu \cite{XiaoXu} generalized the definition of the derived Hall algebra to any triangulated category with some homological finiteness conditions. In \cite{SX}, Sheng and Xu proved that for each hereditary algebra $A$ the lattice algebra of $A$ is isomorphic to its twisted and extended derived Hall algebra.

In 2013, for each finite dimensional hereditary algebra $A$ over a finite field, Bridgeland \cite{Br} defined an algebra associated to $A$, called Bridgeland's Hall algebra of \emph{A}, which is the Ringel--Hall algebra of $2$-cyclic complexes over projective $A$-modules with some localization and reduction. He proved that the quantized enveloping algebra associated to the hereditary algebra \emph{A} can be embedded into Bridgeland's Hall algebra of $A$. This provides a realization of the full quantized enveloping algebra by Hall algebras. Afterwards, Yanagida \cite{Yan} showed that Bridgeland's Hall algebra of $2$-cyclic complexes of any finite dimensional hereditary algebra is isomorphic to the (reduced) Drinfeld double of its extended Ringel--Hall algebra. Inspired by the work of Bridgeland, Chen and Deng \cite{ChenD} considered Bridgeland's Hall algebra $\mathcal {D}\mathcal {H}_m(A)$ of $m$-cyclic complexes of a hereditary algebra $A$ for each non-negative integer $m\neq1$.

In this paper, let $A$ be a finite dimensional hereditary algebra over a finite field. For any non-negative integer $m\neq1$ or $2$ we first define an $m$-periodic lattice algebra, and then use it to give a characterization of the algebra structure on $\mathcal {D}\mathcal {H}_m(A)$. Explicitly, we show that Bridgeland's Hall algebra $\mathcal {D}\mathcal {H}_m(A)$ is isomorphic to the $m$-periodic lattice algebra. As a byproduct, we show that there is an embedding of the Heisenberg double Hall algebra of $A$ into $\mathcal {D}\mathcal {H}_m(A)$.

Throughout the paper, let $m$ be a non-negative integer such that $m=0$ or $m>2$, let $k$ be a fixed finite field with $q$ elements and set $v=\sqrt{q}\in \mathbb{C}$. Denote by $A$ a finite dimensional $k$-algebra. We denote by $\mod A$ and $D^b(A)$ the category of finite dimensional (left) $A$-modules and its bounded derived category, respectively, and denote by $\mathscr{P}=\mathscr{P}_A$ the full subcategory of $\mod A$ consisting of projective $A$-modules. Let $K(A)$ be the Grothendieck group of $\mod A$ and $\Iso( A)$ the set of isoclasses (isomorphism classes) of $A$-modules. For an $A$-module $M$, the class of $M$ in $K(A)$ is denoted by $\hat{M}$, and the automorphism group of $M$ is denoted by $\Aut(M)$. For a finite set $S$, we denote by $|S|$ its cardinality. We also write $a_M$ for $|\Aut(M)|$.
For a complex $M_\bullet$ of $A$-modules, its homology is denoted by $H_\ast(M_\bullet)$. For a positive integer $m$, we denote the quotient ring $\mathbb{Z}/m\mathbb{Z}$ by $\mathbb{Z}_m=\{0,1,\ldots,m-1\}$. By convention, $\mathbb{Z}_0=\mathbb{Z}$.

\section{Preliminaries}

In this section, we summarize some necessary definitions and properties. All of the materials can be found in \cite{Br}, \cite{ChenD}, \cite{Sch} and \cite{ZHC1}.
\subsection{$m$-cyclic complexes}
Let $m$ be a positive integer. Given an additive category $\mathcal{A}$, an \emph{$m$-cyclic complex} ${{M}_{\bullet }}={{({{M}_{i}},{{d}_{i}})}_{i\in {{\mathbb{Z}}_{m}}}}$ over $\mathcal{A}$ consists of $m$ objects $M_i$ and $m$ morphisms $d_i:M_i\rightarrow M_{i+1}$ in $\mathcal{A}$
satisfying $d_{i+1}d_i=0$ for all $i\in \mathbb{Z}_m$. Hence, each $m$-cyclic complex ${{M}_{\bullet }}={{({{M}_{i}},{{d}_{i}})}_{i\in {{\mathbb{Z}}_{m}}}}$ can be diagrammed by $$\xymatrix@!=0.05cm @M=0pt{&&M_0\ar@/^/[dr]^-{d_0}&\\
&M_{m-1}\ar@/^/[ur]^--{d_{m-1}}&&M_1\ar@/^0.5pc/[dl]\\
&&\ar@/^/[ul]\cdots &}$$ with $d_{i+1}d_i=0$ for all $i\in \mathbb{Z}_m$. A \emph{morphism} $f$ between two $m$-cyclic complexes ${{M}_{\bullet }}={{({{M}_{i}},{{d}_{i}})}_{i\in {{\mathbb{Z}}_{m}}}}$ and
${{N}_{\bullet }}={{({{N}_{i}},{{c}_{i}})}_{i\in {{\mathbb{Z}}_{m}}}}$ is given by $m$ morphisms $f_i:M_i\rightarrow N_i$ in $\mathcal{A}$ satisfying $f_{i+1}d_i=c_if_i$
for all $i\in\mathbb{Z}_m$.

Let $f=(f_i)_{i\in\mathbb{Z}_m}$ and $g=(g_i)_{i\in\mathbb{Z}_m}$ be two morphisms between $m$-cyclic complexes ${{M}_{\bullet }}={{({{M}_{i}},{{d}_{i}})}_{i\in {{\mathbb{Z}}_{m}}}}$ and
${{N}_{\bullet }}={{({{N}_{i}},{{c}_{i}})}_{i\in {{\mathbb{Z}}_{m}}}}$.
We say that $f$ is \emph{homotopic} to $g$ if there exist $m$ morphisms $s_i:M_i\rightarrow N_{i-1}$ in $\mathcal{A}$ such that $f_i-g_i=s_{i+1}d_i+c_{i-1}s_i$ for all $i\in\mathbb{Z}_m$. The category of $m$-cyclic complexes over $\mathcal{A}$ is denoted by $\mathcal{C}_m(\mathcal{A})$, and $K_m(\mathcal{A})$ denotes the homotopy category of~$\mathcal{C}_m(\mathcal{A})$ by identifying homotopic morphisms. As in usual complex categories, one can define quasi-isomorphisms in $\mathcal{C}_m(\mathcal{A})$ and $K_m(\mathcal{A})$, and then get a triangulated category, denoted by $\mathcal{D}_m(\mathcal{A})$, by localizing $K_m(\mathcal{A})$ with respect to the set of all quasi-isomorphisms. We write $\mathcal{C}_0(\mathcal{A})$,$K_0(\mathcal{A})$ and $\mathcal{D}_0(\mathcal{A})$ for the category of bounded complexes, its homotopy category and derived category, respectively.

For each integer $t$, there is a shift functor $$[t]:\mathcal{C}_m(\mathcal{A})\rightarrow \mathcal{C}_m(\mathcal{A}),~~M_\bullet \mapsto M_\bullet[t],$$
where $M_\bullet[t]=(X_i,f_i)$ is defined by $$X_i=M_{i+t},~~f_i=(-1)^{t}d_{i+t},~~ i\in {\mathbb{Z}_m}.$$

Recall that $A$ is a finite dimensional $k$-algebra throughout the paper. Applying the above construction to $\mathscr{P}=\mathscr{P}_A$, we obtain the categories $\mathcal{C}_m(\mathscr{P})$ and $K_m(\mathscr{P})$. In the sense of chain-wise exactness, $\mathcal{C}_m(\mathscr{P})$ is a Frobenius exact category.

For an arbitrary homomorphism $f:Q\ra P$~of projective $A$-modules, define $C_f=(X_i,d_i)\in \mathcal{C}_m(\scrP)$ by
$$X_i=\begin{cases} Q\;\;&\text{~$i=m-1$};\\
                    P\;\;&\text{~$i=0$};\\
                     0 &\text{otherwise,}\end{cases}\qquad
 d_i=\begin{cases} f\;\;&\text{~$i=m-1$};\\
                     0 &\text{otherwise.}\end{cases}$$
So each projective $A$-module $P$ determines an object ${{K}_{P}}:=C_{id_P}$ in $\mathcal{C}_m(\mathscr{P})$.

The following lemma taken from \cite{ChenD} is important in the later calculations.
\begin{lemma}{\rm(\cite[Lem. 2.5]{ChenD})}\label{Ext to Hom}
If $M_\bullet,N_\bullet \in \mathcal{C}_m(\mathscr{P})$, then there exists an isomorphism of vector spaces $$\Ext_{\mathcal{C}_m(\mathscr{P})}^1(N_\bullet,M_\bullet) \cong \Hom_{K_m(\mathscr{P})}(N_\bullet,M_\bullet[1]).$$
\end{lemma}

Let $F$ be the natural covering functor from $\mathcal{C}_0(\mod A)$ to $\mathcal{C}_m(\mod A)$ with Galois group $\lr{[m]}$. Let $D^b(A)/[m]$ be the orbit category of $D^b(A)$ with respect to the shift functor $[m]$.

\begin{lemma}
If $A$ is of finite global dimension, then
$F$ induces a fully faithful functor $$F:D^b(A)/[m]\longrightarrow K_m(\mathscr{P}).$$
\end{lemma}
\bp
Since $A$ is of finite global dimension, we can equally well define $D^b(A)/[m]$ as the orbit category of $K^b(\mathscr{P})$. Then it is easy to see that $F$ induces a fully faithful functor $F:K^b(\mathscr{P})/[m]\longrightarrow K_m(\mathscr{P}).$
\ep
An $m$-cyclic complex $M_\bullet$ of $A$-modules is called \emph{acyclic} if $H_\ast(M_\bullet)=0$.
For each $A$-module $P\in\mathscr{P}$, $K_P[r]$ for all $r\in\mathbb{Z}_m$ are acyclic.
Clearly, for any complex $M_\bullet\in \mathcal{C}_m(\mathscr{P})$, $M_\bullet$ is acyclic if $M_\bullet \cong 0$ in $K_m(\mathscr{P})$.
The following lemma gives a characterization of all acyclic complexes in $\mathcal{C}_m(\mathscr{P})$.

\begin{lemma}{\rm(\cite[Lem. 3.2]{Br};\cite[Lem. 2.2]{ZHC1})}
Suppose that $A$ is of finite global dimension. Then for each acyclic complex $M_\bullet\in \mathcal{C}_m(\mathscr{P})$, there exist objects $P_r\in\mathscr{P}$, $r\in\mathbb{Z}_m$, such that ${{M}_{\bullet }}\cong\underset{r\in\mathbb{Z}_m}\bigoplus K_{P_r}[r]$. Moreover, these projective $A$-modules $P_r$ are uniquely determined up to isomorphism.
\end{lemma}

\subsection{Hall algebras}
Given $A$-modules $L, M, N$, we define $\Ext_{A}^1(M,N)_L$ to be the subset of $\Ext_{A}^1(M,N)$, which consists of those equivalence classes of short exact sequences with middle term $L$.
We define the \emph{Hall algebra} $\mathfrak{H}(A)$ to be the vector space over $\mathbb{C}$ with basis $[M]\in \Iso( A)$ and with associative multiplication defined by
\[[M] \diamond [N] = \sum\limits_{[L] \in \Iso( A)} F_{MN}^L[L],\]
where $$F_{MN}^L={\frac{{|\Ext_A^1{{(M,N)}_L}|}}{{|\Hom_A(M,N)|}}}\cdot\frac{|\Aut L|}{|\Aut M|\cdot |\Aut N|}$$ and it is called the \emph{Ringel--Hall number} associated to $A$-modules $L,M,N$.

From now on, we suppose that $A$ is hereditary. For $M, N\in \mod A$, define $$\lr{M,N}:=\dim_k{\Hom_{A}(M,N)}-\dim_k{\Ext_{A}^{1}(M,N)}.$$
It induces a bilinear form
$$\lr{\cdot ,\cdot }: K(A)\times K(A)\longrightarrow \mathbb{Z},$$ known as the \emph{Euler form}. We also consider the \emph{symmetric Euler form}
$$(\cdot ,\cdot ): K(A)\times K(A)\longrightarrow \mathbb{Z},$$ defined by $(\alpha,\beta)=\lr{\alpha,\beta}+\lr{\beta,\alpha}$ for all $\alpha,\beta \in K(A)$.

The twisted Hall algebra $\mathfrak{H}_v(A)$, called the \emph{Ringel--Hall algebra}, is the same vector space as $\mathfrak{H}(A)$ but with twisted multiplication defined by $$[M]\ast[N]=v^{\lr{{M},\,{N}}}\cdot[M]\diamond[N].$$
We can form the \emph{extended Ringel--Hall algebra} $\mathfrak{H}_v^e(A)$ by adjoining symbols $K_\alpha$ for $\alpha\in K(A)$ and imposing relations
\begin{flalign}
K_\alpha\ast K_\beta=K_{\alpha+\beta},~~~~K_\alpha \ast [M]=v^{\lr{\alpha,\,\hat{M}}}\cdot [M] \ast K_\alpha.
\end{flalign}

\begin{definition}
The \emph{twisted Hall algebra} ${\mathcal {H}}_{\rm tw}(\mathcal{C}_m(\mathscr{P}))$ of $\mathcal{C}_m(\mathscr{P})$ is the vector space over $\mathbb{C}$ with basis indexed by the isoclasses $[M_\bullet]$ of objects in ${\mathcal{C}_m(\mathscr{P})}$, and with multiplication defined by
\[[M_\bullet] \ast [N_\bullet] = v^{\lr{{M_\bullet},\,{N_\bullet}}}\cdot\sum\limits_{[L_\bullet]}{\frac{{|\Ext_{\mathcal{C}_m(\mathscr{P})}^1{{(M_\bullet,N_\bullet)}_{L_\bullet}}|}}{{|\Hom_{\mathcal{C}_m(\mathscr{P})}(M_\bullet,N_\bullet)|}}} [L_\bullet],\]where $\lr{{M_\bullet},\,{N_\bullet}}:=\sum\limits_{i\in {{\mathbb{Z}}_{m}}}{\lr{{{\hat{M}}_{i}},\,{{\hat{N}}_{i}}}}$.
\end{definition}
It is easy to see that $\mathcal {H}_{\rm tw}(\mathcal{C}_m(\mathscr{P}))$ is an associative algebra (cf. \cite{Br,ChenD}). By some simple calculations, we have the following relations for the acyclic complexes ${K_P[r]}$, $r\in\mathbb{Z}_m$.

\begin{lemma}{\rm(\cite{ZHC1})}\label{Ore}
Let $M\in \mod A$, $P, Q \in \mathscr{P}$ and ${M_\bullet}\in\mathcal{C}_m(\mathscr{P})$. We have the following identities for each $r\in{{\mathbb{Z}}_{m}}$  in $\mathcal {H}_{\rm tw}(\mathcal{C}_m(\mathscr{P}))$
\begin{flalign}
&[K_P[r]]\ast[M_\bullet]=v^{\lr{\hat{P},\,\hat{M}_{m-r}-\hat{M}_{m-r-1}}}[K_P[r] \oplus M_\bullet];\\
&[M_\bullet]\ast[K_P[r]]=v^{-\lr{\hat{M}_{m-r}-\hat{M}_{m-r-1},\,\hat{P}}}[K_P[r] \oplus M_\bullet];\\
&[K_P[r]]\ast[M_\bullet]=v^{(\hat{P},\,\hat{M}_{m-r}-\hat{M}_{m-r-1})}[M_\bullet]\ast[K_P[r]]\label{jiaohuan1}.
\end{flalign}
In particular,

\begin{equation}\label{jiaohuan3}{[K_P[r]]}\ast {[K_Q]}=\begin{cases}
v^{( \hat{P},\,\hat{Q})}{[K_Q]}\ast {[K_P[r]]} \quad &\text{~$r=1$},\\
v^{-( \hat{P},\,\hat{Q})}{[K_Q]}\ast {[K_P[r]]} \quad &\text{~$r=m-1$},\\
{[K_Q]}\ast {[K_P[r]]} & {otherwise;}
\end{cases}\end{equation}
\begin{equation}\label{jiaohuan4}{[K_P[r]]}\ast {[C_M]}=\begin{cases}
v^{( \hat{P}, \hat{M})}{[C_M]}\ast {[K_P[r]]} \quad &\text{~$r=0$},\\
v^{( \hat{P},\,  \hat{\ooz}_M)}{[C_M]}\ast {[K_P[r]]} \quad &\text{~$r=1$},\\
v^{-( \hat{P},\,  \hat{P}_M)}{[C_M]}\ast {[K_P[r]]} \quad &\text{~$r=m-1$},\\
{[C_M]}\ast {[K_P[r]]} &{otherwise}.\\
\end{cases}\end{equation}
\end{lemma}

Set $\mathscr{S}:=\left\{v^n[{{M}_{\bullet }}]\in {{\mathcal{C}}_{m}}(\mathscr{P})~|~{{ {H}}_{*}}({{M}_{\bullet }})=0,~ n\in\mathbb{Z}\right\}$. By Lemma
\ref{Ore}, the set $\mathscr{S}$ satisfies the Ore conditions. So one considers the localization of $\mathcal {H}_{\rm tw}(\mathcal{C}_m(\mathscr{P}))$ with respect to the set $\mathscr{S}$ (cf. \cite{Br,ChenD}).

\begin{definition}
The localized Hall algebra $\mathcal {D}\mathcal {H}_m(A)$ of $A$, called \emph{Bridgeland's Hall algebra}, is the localization of $\mathcal {H}_{\rm tw}(\mathcal{C}_m(\mathscr{P}))$ with respect to the set $\mathscr{S}$.
Since $v$ is invertible,
in symbols, $$\mathcal {D}\mathcal {H}_m(A):=\mathcal {H}_{\rm tw}(\mathcal{C}_m(\mathscr{P}))[~[M_\bullet]^{-1}~|~H_\ast(M_\bullet)=0~].$$

\end{definition}

For any $\alpha \in K(A)$ and $r\in \mathbb{Z}_m$, by writing $\alpha$ in the form $\alpha =\hat{P}-\hat{Q}$ for some $A$-modules $P,Q\in\mathscr{P}$, we define
$${{K}_{\alpha ,r}}:=[{{K}_{P}}[r]]*{{[K_Q[r]]}^{-1}},$$
and it is easy to see that \begin{equation}\label{kejia}
K_{\alpha,r}\ast K_{\beta,r}=K_{\alpha+\beta,r}.\end{equation}
We simply write $K_\alpha$ for $K_{\alpha,0}$.
Note that the identities $(\ref{jiaohuan1})$--$(\ref{jiaohuan4})$ continue to hold with the elements $[K_P[r]]$ and $[K_Q]$ replaced by $K_{\alpha,r}$ and $K_\beta$, respectively, for all $\alpha, \beta \in K(A)$.

Each $A$-module $M$ has a unique minimal projective resolution up to isomorphism of the form \footnote{For each fixed $A$-module $M$, we fix a minimal projective resolution (\ref{projective resolution}) of $M$ using notations $P_M$ and $\Omega_M$ throughout the paper.}
\begin{equation}\label{projective resolution}
0\lra \ooz_M\stackrel{\delta_M}{\lra} P_M\stackrel{\vez_M}{\lra}
M\lra 0.
\end{equation}
Given an $A$-module $M$,
we take a minimal projective resolution (\ref{projective resolution}) of $M$, and consider the corresponding $m$-cyclic complex
$C_M:=C_{\delta_M}$.
Since the uniqueness of the minimal projective resolution up to isomorphism, the complex $C_M$ is well-defined up to isomorphism.
By \cite{Br,ChenD}, for each $r\in\mathbb{Z}_m$, we define  $$E_{M,r} := v^{\lr{\hat{\Omega}_M,\hat{M}}} \cdot {K_{-\hat{\Omega}_{M,r}}} \ast [C_M[r]] \in \mathcal {D}\mathcal {H}_m(A).$$

For each $r\in\mathbb{Z}_m$, set $e_{M,r}:=a_M^{-1}\cdot E_{M,r}$. We also simply write $E_M$ and $e_M$ for $E_{M,0}$ and $e_{M,0}$, respectively. Let us reformulate a result from \cite[Prop. 4.4]{ChenD}.
\begin{proposition}\label{CD}
There is an embedding of algebras for each $i\in\mathbb{Z}_m$
\begin{equation}\label{emb} I_i:\mathfrak{H}_{v}^e(A)\hookrightarrow\mathcal {D}\mathcal {H}_m(A),~~[M]\mapsto e_{M,i},~~K_\alpha\mapsto K_{\alpha,i}. \end{equation} Moreover, the multiplication map induces an isomorphism of vector spaces
\begin{equation}\label{ji}\begin{aligned}
&\bigotimes\limits_{i\in\mathbb{Z}_m}
\mathfrak{H}_{v}^e(A)\longrightarrow \mathcal {D}\mathcal {H}_m(A),\quad
&\bigotimes\limits_{i\in\mathbb{Z}_m}x_i\longmapsto \prod\limits_{i\in\mathbb{Z}_m}^{\rightarrow}I_i(x_i).
\end{aligned}\end{equation}
\end{proposition}

\begin{remark}
If $m=0$, then the tensor product in Proposition \ref{CD} is an infinite tensor. Here and elsewhere in this paper all infinite tensor products are understood in the restricted sense (cf. \cite[Sec. 3.2]{Kapranov}). $\prod\limits_{i\in\mathbb{Z}_m}^{\rightarrow}a_i$ means the ordered product of the elements $a_i$: $\prod\limits_{i\in\mathbb{Z}_m}^{\rightarrow}a_i=
\cdots a_i\cdot a_{i+1}\cdot\cdots$.
\end{remark}

\section{Heisenberg doubles}
Recall that $A$ is a finite dimensional hereditary $k$-algebra. First of all, we give a counting symbol.
For any fixed $M,N,X,Y\in\mod A$, we denote by $W_{MN}^{XY}$ the set $$\{(f,g,h)~|~\xymatrix{0\ar[r]&X\ar[r]^g&M\ar[r]^f&N\ar[r]^h&Y\ar[r]&0} \text{is exact in}~~\mod A\},$$ and set $$\gamma_{MN}^{XY}:=\frac{|W_{MN}^{XY}|}{a_Ma_N}.$$
Let $\Heis(A)$ be the \emph{Heisenberg double} of the extended Ringel--Hall algebra $\mathfrak{H}_{v}^e(A)$ (cf. \cite[Sec. 1.5]{Kapranov}). By definition, $\Heis(A)$ is an associative and unital $\mathbb{C}$-algebra generated by elements $Z_{[M]}^+, Z_{[M]}^{-}$ and $K_\alpha, K_\alpha^{-}$ with
$[M]\in\Iso(A)$ and $\alpha\in K(A)$, which are subject to the following relations:
\begin{equation*}\begin{split}
&Z_{[M]}^{+}Z_{[N]}^{+}=v^{\lr{M,\,N}}\sum\limits_{[L]}F_{MN}^LZ_{[L]}^{+},~~
K_\alpha Z_{[M]}^{+}=v^{(\alpha,\,\hat{M})}Z_{[M]}^{+}K_\alpha;\\
&Z_{[M]}^{-}Z_{[N]}^{-}=v^{\lr{M,\,N}}\sum\limits_{[L]}F_{MN}^LZ_{[L]}^{-},~~
K_\alpha^- Z_{[M]}^{-}=v^{(\alpha,\,\hat{M})}Z_{[M]}^{-}K_\alpha^-;\\
&Z_{[M]}^{-}Z_{[N]}^{+}=\sum\limits_{[X],[Y]}v^{\lr{\hat{M}-\hat{X},\,\hat{X}-\hat{Y}}}\gamma_{MN}^{XY}K_{\hat{M}-\hat{X}}Z_{[Y]}^{+}Z_{[X]}^{-};\\
&K_\alpha Z_{[M]}^{-}=v^{-(\alpha,\,\hat{M})}Z_{[M]}^{-}K_\alpha,~~K_\alpha^-Z_{[M]}^{+}=Z_{[M]}^{+}K_\alpha^-,~~K_\alpha K_\beta^-=v^{-(\alpha,\,\beta)}K_\beta^-K_\alpha.
\end{split}\end{equation*}
Clearly, $\Heis(A)$ is naturally related to the subcategory $D^{[-1,0]}(A)$ of the derived category $D^b(A)$, which is consisting of two copies of $\mod A$ inside $D^b(A)$ given by complexes concentrated in degrees $0$ and $-1$. In other words, $\Heis(A)$ gives rise to  two copies of $\mathfrak{H}_{v}^e(A)$ with Heisenberg double-type commutation relations (cf. \cite[Prop. 1.5.3]{Kapranov}). Moreover, Kapranov \cite[Def. 3.1]{Kapranov} introduced an associative and unital $\mathbb{C}$-algebra $L(A)$, called the \emph{lattice algebra} of $A$, by taking not just two but infinitely many copies of $\mathfrak{H}_{v}(A)$ and one copy of the group algebra $\mathbb{C}[K(A)]$, and by imposing Heisenberg double-like commutation relations between adjacent copies of $\mathfrak{H}_{v}(A)$ and oscillator relations of the form $ab=\lambda_{ab}ba,~\lambda_{ab}\in\mathbb{R}$, between basis vectors of non-adjacent copies. In a similar manner to \cite[Def. 3.1]{Kapranov}, we give the following definition.
\begin{definition}\label{m-periodic}
The \emph{$m$-periodic lattice algebra} $L_m(A)$ of $A$ is the associative and unital $\mathbb{C}$-algebra generated by the elements in $\{Z_{[M]}^{(i)}~|~[M]\in\Iso( A),~i\in \mathbb{Z}_m\}$ and $\{K_\alpha^{(i)}~|~\alpha\in K(A),~i\in \mathbb{Z}_m\}$ with the following relations:
\begin{equation}\label{x1}
{K_\alpha^{(i)}} {K_\beta^{(i)}}={K_{\alpha+\beta}^{(i)}},~~~~
{K_\alpha^{(i)}} {K_\beta^{(j)}}=\begin{cases}
v^{( \alpha,\,\beta)}\cdot{K_\beta^{(j)}} {K_\alpha^{(i)}} \quad &\text{~$i=j+1$},\\
v^{-( \alpha,\,\beta)}\cdot{K_\beta^{(j)}} {K_\alpha^{(i)}} \quad &\text{~$i=m-1+j$},\\
{K_\beta^{(j)}} {K_\alpha^{(i)}} & {otherwise;}
\end{cases}\end{equation}
\begin{equation}\label{x2}
{K_\alpha^{(i)}}Z_{[M]}^{(j)}=\begin{cases}
v^{(\alpha,\,\hat{M})}\cdot Z_{[M]}^{(j)}K_\alpha^{(i)} \quad &\text{~$i=j$},\\
v^{-(\alpha,\,\hat{M})}\cdot Z_{[M]}^{(j)}K_\alpha^{(i)} \quad &\text{~$i=m-1+j$},\\
Z_{[M]}^{(j)}K_\alpha^{(i)} & {otherwise;}
\end{cases}\end{equation}
\begin{equation}\label{x3}
Z_{[M]}^{(i)}Z_{[N]}^{(i)}=v^{\lr{M,\,N}}\sum\limits_{[L]}F_{MN}^LZ_{[L]}^{(i)};
\end{equation}
\begin{equation}\label{x4}
Z_{[M]}^{(i+1)}Z_{[N]}^{(i)}=\sum\limits_{[X],[Y]}v^{\lr{\hat{M}-\hat{X},\,\hat{X}-\hat{Y}}}\gamma_{MN}^{XY}K_{\hat{M}-\hat{X}}^{(i)}Z_{[Y]}^{(i)}Z_{[X]}^{(i+1)};
\end{equation}
\begin{equation}\label{x5}
Z_{[M]}^{(i)}Z_{[N]}^{(j)}=Z_{[N]}^{(j)}Z_{[M]}^{(i)},~~i-j\neq0, 1~\mbox{or}~m-1.
\end{equation}
\end{definition}

\begin{theorem}\label{main}
Bridgeland's Hall algebra $\mathcal {D}\mathcal {H}_m(A)$ is isomorphic to the $m$-periodic lattice algebra $L_m(A)$.
\end{theorem}
\bp
First of all, we claim that there is a homomorphism of algebras
$$\Phi: L_m(A)\longrightarrow \mathcal {D}\mathcal {H}_m(A)$$
defined by $\Phi(K_\alpha^{(i)})=K_{\alpha,i}$ and $\Phi(Z_{[M]}^{(i)})=e_{M,i}$. That is to say, $\Phi$ preserves the relations $(\ref{x1})$--$(\ref{x5})$. Indeed, by the identities (\ref{kejia}) and $(\ref{jiaohuan3})$, clearly, the relation $(\ref{x1})$ is preserved.
By the identities $(\ref{jiaohuan3})$ and $(\ref{jiaohuan4})$, for any $M\in \mod A$ and $P\in\mathscr{P}$, it is easy to see that
\begin{equation}
[K_P[r]]\ast (K_{-\hat{\Omega}_M}\ast[C_M])=\begin{cases}
v^{(P,\,M)}\cdot (K_{-\hat{\Omega}_M}\ast[C_M])\ast [K_P[r]]\quad &\text{~$r=0$},\\
v^{-(P,\,M)}\cdot (K_{-\hat{\Omega}_M}\ast[C_M])\ast [K_P[r]]\quad &\text{~$r=m-1$},\\
(K_{-\hat{\Omega}_M}\ast[C_M])\ast [K_P[r]] & {otherwise.}
\end{cases}\end{equation}
Hence, the relation $(\ref{x2})$ is preserved. By Proposition \ref{CD}, the relation $(\ref{x3})$ is preserved.
Let us first consider the relation (\ref{x5}). Let $M,N\in\mod A$ and $r\neq 0$ or 1 in $\mathbb{Z}_m$, then
\begin{equation}\begin{split}
\Ext^1_{\mathcal{C}_m(\mathscr{P})}(C_M[r],C_N)&\cong\Hom_{K_m(\mathscr{P})}(C_M[r],C_N[1])\\
&\cong\Hom_{K_m(\mathscr{P})}(C_M[r-1],C_N)\\&=0.\end{split}\end{equation}
Similarly, for $r\neq 0$ or $m-1$ in $\mathbb{Z}_m$, we obtain that $\Ext^1_{\mathcal{C}_m(\mathscr{P})}(C_N,C_M[r])=0$. So it is easy to see that for $r\neq0, 1$ or $m-1$, we have $$[C_M[r]]\ast[C_N]=[C_N]\ast[C_M[r]]=[C_M[r]\oplus C_N].$$
Hence, \begin{equation}\begin{split}
(K_{-\hat{\Omega}_{M,r}}\ast[C_M[r]])\ast(K_{-\hat{\Omega}_N}\ast[C_N])&=
K_{-\hat{\Omega}_{M,r}}\ast K_{-\hat{\Omega}_N}\ast[C_M[r]]\ast[C_N]\\
&=K_{-\hat{\Omega}_{M,r}}\ast K_{-\hat{\Omega}_N}\ast[C_N]\ast[C_M[r]]\\
&=K_{-\hat{\Omega}_N}\ast K_{-\hat{\Omega}_{M,r}}\ast[C_N]\ast[C_M[r]]\\
&=(K_{-\hat{\Omega}_N}\ast[C_N])\ast (K_{-\hat{\Omega}_{M,r}}\ast[C_M[r]]).
\end{split}\end{equation}
Namely, $e_{M,r}\ast e_N=e_N\ast e_{M,r}$. So the relation $(\ref{x5})$ is preserved. We remove the proof of the relation $(\ref{x4})$ to the next section and prove that $\Phi$ is an isomorphism right now.

For each $i\in\mathbb{Z}_m$, let $L_{m}^{(i)}(A)$ be the subalgebra of $L_{m}(A)$ generated by the elements in  $\{Z_{[M]}^{(i)}~|~[M]\in\Iso( A)\}$ and $\{K_\alpha^{(i)}~|~\alpha\in K(A)\}$. Then, clearly, the multiplication map of $L_{m}(A)$ induces an epimorphism of vector spaces
$$\xymatrix{\Theta:\bigotimes\limits_{i\in\mathbb{Z}_m}L_m^{(i)}(A)\ar@{->>}[r]& L_{m}(A).}$$

For each $i\in\mathbb{Z}_m$, let $\mathcal{H}^{(i)}(A)$ be the subalgebra of $\mathcal {D}\mathcal {H}_m(A)$ generated by the elements in  $\{e_{M,i}~|~[M]\in\Iso( A)\}$ and $\{K_{\alpha,i}~|~\alpha\in K(A)\}$.
Then, by Proposition \ref{CD}, for each $i\in\mathbb{Z}_m$, $\Phi$ induces an algebra isomorphism $$\Phi_i: L_{m}^{(i)}(A)\longrightarrow \mathcal{H}^{(i)}(A),$$ moreover, there exists an isomorphism of vector spaces
$$\Psi : \xymatrix{\bigotimes\limits_{i\in\mathbb{Z}_m}\mathcal{H}^{(i)}(A)\ar[r]&\mathcal {D}\mathcal {H}_m(A).}$$

Hence, we have the following commutative diagram
$$\xymatrix{\bigotimes\limits_{i\in\mathbb{Z}_m}L_m^{(i)}(A)\ar@{->>}[r]^-{ \Theta}\ar[d]_-{\bigotimes\limits_{i\in\mathbb{Z}_m}\Phi_{i}}^-{\cong} &L_m(A)\ar[d]^-{\Phi}\\
\bigotimes\limits_{i\in\mathbb{Z}_m}\mathcal{H}^{(i)}(A)\ar[r]^-{\ \cong} &\mathcal {D}\mathcal {H}_m(A).}$$
Since $\Phi\circ\Theta$ is an isomorphism, we obtain that $\Theta$ is injective, and thus it is an isomorphism. So $\Phi$ is an isomorphism.
\ep

\begin{corollary}
For each $i\in\mathbb{Z}_m$, there exists an embedding of algebras
\begin{equation*}\label{emb2} J_i:\Heis(A)\hookrightarrow\mathcal {D}\mathcal {H}_m(A),~~Z_{[M]}^{+}\mapsto e_{M,i},~~Z_{[M]}^{-}\mapsto e_{M,i+1},~~K_\alpha\mapsto K_{\alpha,i},~~K_\alpha^{-}\mapsto K_{\alpha,i+1}. \end{equation*}
\end{corollary}
\bp
For each $i\in\mathbb{Z}_m$, clearly, the map
\begin{equation*} J'_i:\Heis(A)\hookrightarrow L_m(A),~~Z_{[M]}^{+}\mapsto Z_{[M]}^{(i)},~~Z_{[M]}^{-}\mapsto Z_{[M]}^{(i+1)},~~K_\alpha\mapsto K_{\alpha}^{(i)},~~K_\alpha^{-}\mapsto K_{\alpha}^{(i+1)}\end{equation*} is an embedding of algebras.
\ep

\section{The proof of the relation (\ref{x4}) in Theorem \ref{main}}
Consider an extension of $C_M[1]$ by $C_N$
$$\eta: 0\lra C_N\lra L_{\bullet}\lra C_M[1] \lra 0.$$
It induces a long exact sequence in homology
$$H_{m-1}(C_N)\lra H_{m-1}(L_{\bullet})\lra H_{m-1}(C_M[1])\lra H_0(C_N)\lra H_0(L_{\bullet})\lra H_0(C_M[1]).$$ Clearly, $H_{m-1}(C_Z)=H_0(C_Z[1])=0$ and $H_{m-1}(C_Z[1])=H_0(C_Z)=Z$ for any $Z\in\mod A$. Hence, by writing
$$L_\bullet=C_X[1]\oplus C_Y\oplus K_T\oplus K_W[1]$$ for some $X, Y\in
\mod A$ and $T, W\in\scrP$,
we obtain an exact sequence of $A$-modules
\begin{equation}\label{long-exact-sequence}
0\lra X\lra M\stackrel{\delta}{\lra} N\lra Y\lra 0,
\end{equation}where $\delta$ is determined by the equivalence class of $\eta$ via the canonical isomorphisms
$$\Ext^1_{\mathcal{C}_m(\scrP)}(C_M[1],C_N)\cong \Hom_{K_m(\scrP)}(C_M[1],C_N[1])\cong \Hom_A(M,N).$$ Clearly, $\eta $ splits if and only if $\delta =0$.
By considering the kernels and cokernels of differentials in $C_N, C_M[1]$ and $L_{\bullet}$, we obtain that \begin{equation}\label{TW}T\oplus P_Y\cong P_N,~W\oplus \ooz_X\cong \ooz_M.\end{equation}
It is easy to see that for any $X, Y\in
\mod A$ and $T, W\in\scrP$, \begin{equation}|\Ext^1_{\mathcal{C}_m(\scrP)}(C_M[1],C_N)_{C_X[1]\oplus C_Y\oplus K_T\oplus K_W[1]}|=\frac{|W_{MN}^{XY}|}{a_Xa_Y}=\frac{a_Ma_N}{a_Xa_Y}\gamma_{MN}^{XY}.\end{equation}

\noindent $\mathbf{Proof~of~Theorem~\ref{main}}$ Let $M,N\in\mod A$.
\begin{equation*}\begin{split}&[C_M[1]]\ast[C_N]=\\&v^{\lr{P_M,\,\Omega_N}}\sum\limits_{[X],[Y]}
|\Ext^1_{\mathcal{C}_m(\scrP)}(C_M[1],C_N)_{C_X[1]\oplus C_Y\oplus K_T\oplus K_W[1]}|\cdot[C_X[1]\oplus C_Y\oplus K_T\oplus K_W[1]].\end{split}\end{equation*}
It is easy to see that $[K_T\oplus K_W[1]]\ast[C_X[1]\oplus C_Y]=v^{x_0-2x_1}\cdot[C_X[1]\oplus C_Y\oplus K_T\oplus K_W[1]]$, where $x_0=\lr{\hat{W},\,\hat{P}_X+\hat{\Omega}_X+\hat{\Omega}_Y}+\lr{\hat{T},\,\hat{P}_X+\hat{P}_Y+\hat{\Omega}_Y}$ and $x_1=\lr{\hat{T},\,\hat{P}_X+\hat{\Omega}_Y}+\lr{\hat{W},\,\hat{\Omega}_X}$. So we obtain that
$$[C_X[1]\oplus C_Y\oplus K_T\oplus K_W[1]]=v^{2x_1-x_0}\cdot[K_T\oplus K_W[1]]\ast[C_X[1]\oplus C_Y].$$
Clearly, $[K_W[1]]\ast[K_T]=v^{\lr{W,\,T}}\cdot[K_W[1]\oplus K_T]$, thus we get that
$$[K_W[1]\oplus K_T]=v^{-\lr{W,\,T}}\cdot[K_W[1]]\ast[K_T].$$
Since $\Ext_{\mathcal{C}_m(\scrP)}^1(C_Y,C_X[1])\cong\Hom_{K_m(\scrP)}(C_Y,C_X[2])=0$,
we have
\begin{equation*}\begin{split}[C_Y]\ast[C_X[1]]&=\frac{v^{\lr{\Omega_Y,\,P_X}}}{|\Hom_{\mathcal{C}_m(\scrP)}(C_Y,C_X[1])|}\cdot[C_X[1]\oplus C_Y]\\&=
v^{-\lr{\Omega_Y,\,P_X}}\cdot[C_X[1]\oplus C_Y].\end{split}\end{equation*}
Hence, $[C_X[1]\oplus C_Y]=v^{\lr{\Omega_Y,\,P_X}}\cdot[C_Y]\ast[C_X[1]].$
Therefore, $$[C_M[1]]\ast[C_N]=v^{x_2}\sum\limits_{[X],[Y]}\frac{a_Ma_N}{{a_X}a_Y}\gamma_{MN}^{XY}
\cdot[K_W[1]]\ast[K_T]\ast[C_Y]\ast[C_X[1]],$$ where $x_2=\lr{P_M,\,\Omega_N}+2x_1-x_0-\lr{W,\,T}+\lr{\Omega_Y,\,P_X}$. So,
\begin{equation*}\begin{split}
e_{M,1}e_N&=(v^{\lr{\Omega_M,\,M}}a_M^{-1}\cdot K_{-\hat{\Omega}_M,1}\ast[C_M[1]])\ast
(v^{\lr{\Omega_N,\,N}}a_N^{-1}\cdot K_{-\hat{\Omega}_N}\ast[C_N])\\
&=v^{\lr{\Omega_M,\,M}+\lr{\Omega_N,\,N}-(P_M,\,\Omega_N)}a_M^{-1}a_N^{-1}\cdot K_{-\hat{\Omega}_M,1}\ast
K_{-\hat{\Omega}_N}\ast[C_M[1]]\ast[C_N]\\&=
v^{x_3}\sum\limits_{[X],[Y]}\frac{\gamma_{MN}^{XY}}{a_Xa_Y}\cdot K_{-\hat{\Omega}_M,1}\ast K_{-\hat{\Omega}_N}
\ast[K_W[1]]\ast[K_T]\ast[C_Y]\ast[C_X[1]]\\&=
v^{x_4}\sum\limits_{[X],[Y]}\frac{\gamma_{MN}^{XY}}{a_Xa_Y}\cdot K_{\hat{T}-\hat{\Omega}_N}
\ast[C_Y]\ast K_{\hat{W}-\hat{\Omega}_M,1}\ast[C_X[1]]
\end{split}\end{equation*}
where
$x_3=\lr{\Omega_M,\,M}+\lr{\Omega_N,\,N}-(P_M,\,\Omega_N)+x_2$,
$x_4=x_3+(\Omega_M,\,\Omega_N)+(\hat{W}-\hat{\Omega}_M,\,\hat{T})+(\hat{W}-\hat{\Omega}_M,\,\hat{\Omega}_Y)$.
Since $\hat{T}-\hat{\Omega}_N=-\hat{\Omega}_Y+\hat{M}-\hat{X}$ and $\hat{W}-\hat{\Omega}_M=-\hat{\Omega}_X$, we obtain that
$$e_{M,1}e_N=\sum\limits_{[X],[Y]}v^{x_5}\gamma_{MN}^{XY}\cdot
K_{\hat{M}-\hat{X}}\ast e_Y\ast e_{X,1},$$ where
$x_5=x_4-\lr{\hat{\Omega}_Y,\,\hat{Y}}-\lr{\hat{\Omega}_X,\,\hat{X}}$.

Using the exact sequence (\ref{long-exact-sequence}), isomorphisms in (\ref{TW}) and the respective minimal projective resolutions (\ref{projective resolution}) of $M,N,X,Y$, we obtain that
\begin{equation*}\begin{split}
x_5&=\lr{\Omega_M,\,M}+\lr{\Omega_N,\,N}-\lr{P_M,\,\Omega_N}-\lr{\Omega_N,\,P_M}+\lr{P_M,\,\Omega_N}
+2\lr{\hat{P}_N-\hat{P}_Y,\,\hat{P}_X+\hat{\Omega}_Y}\\&+
2\lr{\hat{\Omega}_M-\hat{\Omega}_X,\,\hat{\Omega}_X}-\lr{\hat{\Omega}_M-\hat{\Omega}_X,\,\hat{P}_X+\hat{\Omega}_X+\hat{\Omega}_Y}
-\lr{\hat{P}_N-\hat{P}_Y,\,\hat{P}_X+\hat{P}_Y+\hat{\Omega}_Y}\\&
-\lr{\hat{\Omega}_M-\hat{\Omega}_X,\,\hat{P}_N-\hat{P}_Y}+\lr{{\Omega}_Y,\,{P}_X}+\lr{{\Omega}_M,\,{\Omega}_N}+\lr{{\Omega}_N,\,{\Omega}_M}-\lr{{\Omega}_X,\,{\Omega}_Y}-\lr{{\Omega}_Y,\,{\Omega}_X}\\
&-\lr{\hat{\Omega}_X,\,\hat{P}_N-\hat{P}_Y}-\lr{\hat{P}_N-\hat{P}_Y,\,\hat{\Omega}_X}-\lr{{\Omega}_Y,\,Y}-\lr{{\Omega}_X,\,X}
\end{split}
\end{equation*}
\vspace{-0.6em}
\begin{equation*}\begin{split}
&=\lr{\hat{\Omega}_M,\,\hat{M}+2\hat{\Omega}_X-\hat{\Omega}_X-\hat{P}_X-\hat{\Omega}_Y-\hat{P}_N+\hat{P}_Y+\hat{\Omega}_N}
+\lr{\hat{\Omega}_N,\,\hat{N}-\hat{P}_M+\hat{\Omega}_M}\\
&+\lr{\hat{P}_N,\,2\hat{P}_X+2\hat{\Omega}_Y-\hat{P}_X-\hat{\Omega}_Y-\hat{P}_Y-\hat{\Omega}_X}+\lr{\hat{P}_Y,\,-2\hat{P}_X-2\hat{\Omega}_Y+\hat{P}_X+\hat{\Omega}_Y+\hat{P}_Y+\hat{\Omega}_X}
\\&+\lr{\hat{\Omega}_X,\,-2\hat{\Omega}_X+\hat{\Omega}_Y+\hat{P}_X+\hat{\Omega}_X-\hat{P}_Y+\hat{P}_N-\hat{\Omega}_Y-\hat{P}_N-\hat{X}+\hat{P}_Y}+\lr{\hat{\Omega}_Y,\,\hat{P}_X-\hat{\Omega}_X-\hat{Y}}\\
&=\lr{\hat{\Omega}_M,\,\hat{M}-\hat{N}-\hat{X}+\hat{Y}}+\lr{\hat{\Omega}_N,\,\hat{N}-\hat{M}}+\lr{\hat{P}_N,\,\hat{X}-\hat{Y}}+\lr{\hat{P}_Y,\,\hat{Y}-\hat{X}}+\lr{\hat{\Omega}_Y,\,\hat{X}-\hat{Y}}\\
&=\lr{\hat{\Omega}_M-\hat{\Omega}_N,\,\hat{M}-\hat{N}}+\lr{\hat{\Omega}_Y-\hat{P}_Y+\hat{P}_N-\hat{\Omega}_M,\,\hat{X}-\hat{Y}}\\
&=\lr{\hat{\Omega}_M-\hat{\Omega}_N-\hat{Y}+\hat{P}_N-\hat{\Omega}_M,\,\hat{M}-\hat{N}}
=\lr{\hat{N}-\hat{Y},\,\hat{M}-\hat{N}}
=\lr{\hat{M}-\hat{X},\,\hat{X}-\hat{Y}}.
\end{split}
\end{equation*}

Hence, for any $i\in\mathbb{Z}_m$,
$$e_{M,i+1}e_{N,i}=\sum\limits_{[X],[Y]}v^{\lr{\hat{M}-\hat{X},\,\hat{X}-\hat{Y}}}\gamma_{MN}^{XY}\cdot
K_{\hat{M}-\hat{X},i}\ast e_{Y,i}\ast e_{X,i+1}.$$ Therefore, we complete the proof.

\section*{Acknowledgments}

The author is grateful to Bangming Deng and Jie Sheng for their stimulating discussions and valuable comments. He also would like to thank the anonymous referee for modification suggestions.

\end{document}